\documentclass[11pt,a4paper]{article}
\usepackage[utf8]{inputenc}
\usepackage[T1]{fontenc}
\usepackage{lmodern}
\usepackage{amsmath,amssymb,amsthm,mathtools,mathrsfs}
\usepackage{enumitem}
\usepackage{geometry}
\usepackage{hyperref}
\numberwithin{equation}{section}
\hypersetup{colorlinks=true,linkcolor=red,citecolor=blue,urlcolor=blue}

\newcommand{\R}{\mathbb R}

\newcommand{\Cplus}{\mathbb C_+}
\newcommand{\Pp}{\mathcal P}
\newcommand{\Psym}{\mathcal P_{\mathrm{sym}}}
\newcommand{\Pzeroone}{\mathcal P_{0,1}}
\newcommand{\law}{\mathcal L}

\newcommand{\D}{D}

\theoremstyle{plain}
\newtheorem{theorem}{Theorem}[section]

\newtheorem{lemma}[theorem]{Lemma}
\newtheorem{corollary}[theorem]{Corollary}

\theoremstyle{remark}
\newtheorem{remark}[theorem]{Remark}
\newtheorem{example}[theorem]{Example}

\title{Entropic Inequalities and Limit Theorems in Boolean Probability}
\author{Kewei Pan}
\date{}

\begin{document}
\maketitle

\begin{abstract}

We study Boolean entropy and Fisher information in information inequalities and the
Boolean central limit theorem. We prove a fractional Stam inequality for the Boolean Fisher information associated with arbitrary fractional partitions of subset sums, and derive an entropy power inequality for sums of Boolean independent random variables and its weighted version. We then establish a Boolean analogue of Barron's entropy convergence result along the Boolean central limit theorem: Boolean entropy converges to the entropy of the Rademacher law if and only if the initial entropy is finite. Furthermore, under a finite fourth moment, in the Boolean central limit theorem, we obtain
\[
\Gamma(\mu_n\mid\mathbf b)=
\frac{m_4(\mu)-1-m_3(\mu)^2}{2}\frac{\log n}{n}
+o\!\left(\frac{\log n}{n}\right),
\]
with an explicit example showing sharpness. For symmetric laws with finite fourth moment and finite Boolean Fisher information,
\[
\Psi(\mu_n\mid\mathbf b) =2\bigl(m_4(\mu)-1\bigr)\frac{\log n}{n}
+o\left(\frac{\log n}{n}\right).
\]
Our approach is based on the Nevanlinna representation of Cauchy transforms.
\end{abstract}

\tableofcontents

\section{Introduction}
Shannon's entropy of a real-valued random variable $X$ with density $f$ is defined as
\[
 H(X)=-\int_{\R} f(x)\log f(x)\,dx,
\]
whenever the integral is finite. This quantity naturally appears in large deviation theory \cite{sanov1961probability}. Let $X_1,\ldots,X_n$ be a sequence of independent real-valued random variables with finite variance. The classical entropy power inequality of Shannon \cite{shannon1948} and Stam \cite{Stam1959} states that,
\begin{equation}\label{eq:entropy power}
    \exp(2H(X_1+\cdots+X_n))\ge \sum_{i=1}^n \exp(2H(X_i)).
\end{equation}
In \cite{ABBN2004}, the authors showed a stronger version of such an inequality:
\begin{equation}\label{eq:leave one out entropy power}
    \exp(2H(X_1+\cdots+X_n))\ge \frac{1}{n-1}\sum_{i=1}^n \exp(2H(\sum_{j\ne i}X_j)),
\end{equation}
Indeed, taking \(n=2\) in the above inequality yields the entropy power inequality for two summands, which can then be iterated to obtain the corresponding inequality for \(n\) summands. When $X_1,\ldots,X_n$ are identically distributed, using the scaling property: $H(cX)=H(X)+\log |c|$ for any $c\ne0$, \eqref{eq:leave one out entropy power} also implies the monotonicity of Shannon's entropy in the central limit theorem: the entropy of the normalized \(n\)-summand is no smaller than that of the normalized \((n-1)\)-summand. Later, Madiman and Barron \cite{MadimanBarron2007} generalized the entropy power inequality for subset sums. More precisely, if $\mathcal{C}$ is an arbitrary collection of subsets of $\{1,\ldots,n\}$, then
\begin{equation}\label{eq:fractional entropy power}
    \exp(2H(X_1+\cdots+X_n))\ge \frac{1}{r}\sum_{A\in \mathcal{C}} \exp(2H(\sum_{i\in A}X_i)),
\end{equation}
where $r$ is the maximum number of subsets containing any one of the indices. The inequalities \eqref{eq:entropy power} and \eqref{eq:leave one out entropy power} correspond respectively to the families of singletons and of \((n-1)\)-element subsets.

Returning to the i.i.d. setting, suppose that \(X_1\) is centered with variance \(1\). It is well known that the normalized sums converge weakly to the standard Gaussian law \(\gamma\). Moreover, among probability measures with unit second moment, Shannon's entropy \(H(\mu)\) is maximized by \(\gamma\). Recall that the relative entropy with respect to the Gaussian law is defined by
$$
H(\mu\mid\gamma)=\int_\mathbb{R} \log\frac{d\mu}{d\gamma}\,d\mu,
$$
whenever \(\mu\ll\gamma\), and is also known as the Kullback--Leibler divergence. For probability measures with unit second moment, the relative entropy simplifies to
$$
H(\mu\mid\gamma)=H(\gamma)-H(\mu)\geq 0.
$$
The famous Pinsker–Csiszár–Kullback inequality bounds the squared total variation distance to the Gaussian law by the relative entropy. A natural question is therefore whether the convergence of normalized sums to the Gaussian law can be strengthened to convergence at the level of entropy. This question was answered by Barron \cite{Barron1986}, who showed that, if the entropy of the normalized sum is finite for some $n$, then the relative entropy of the normalized sums with respect to the Gaussian law converges to zero. In \cite{ABBN2004rate,JohnsonBarron2004}, quantitative bounds of order $O(n^{-1})$ for the entropy convergence were obtained under the assumption that $X_1$ satisfies a Poincaré inequality. Combined with the Pinsker–Csiszár–Kullback inequality, these estimates yield an $O(n^{-1/2})$ bound for the total variation distance to the Gaussian along the central limit theorem, which is known to be optimal in general. This was further refined by Bobkov, Chistyakov, and Götze \cite{BCG2013rate}, who obtained the sharp constant under the assumption of a finite fourth moment. Their approach is based on Edgeworth-type expansions. These results provide a useful benchmark for investigating analogous questions in other frameworks of independence.

Free probability, developed by Voiculescu \cite{Voiculescu1985symmetries}, provides a successful framework for noncommutative probability, in which the random variables are elements of a unital algebra with a state (unit-preserving linear functional) and freeness plays the role of independence. The distribution of a random variable is encoded by its moments with respect to the state. It is more convenient to work with self-adjoint variables in the context of a $C^*$ or $W^*$-algebra, and their distributions are probability measures on $\mathbb{R}$. The addition of classically independent random variables corresponds to the convolution of probability measures, which is characterized by the additivity of the logarithm of the Fourier transform. For the addition of bounded free random variables, Voiculescu introduced the $R$-transform, which linearizes free additive convolution, and established the free central limit theorem for bounded self-adjoint random variables, with the semicircle law as the limiting distribution \cite{voiculescu1986addition}. The unbounded case was treated separately by Maassen \cite{MAASSEN1992addition} and Bercovici and Voiculescu \cite{Bercovici1993free}. A remarkable discovery of Voiculescu \cite{voiculescu1991limit} is that large unitary invariant random matrices are asymptotically free, and by looking at Shannon's entropy on large random matrices, Voiculescu subsequently introduced free entropy as a free analogue of classical entropy \cite{Voiculescu1993}. For a bounded and self-adjoint random variable $X$ with distribution $\mu$, it is defined as
$$
\Sigma(X)=\iint_{\mathbb{R}^2} \log|x-y|\,d\mu(x)d\mu(y)\in [-\infty,\infty).
$$
In \cite{szarek1996volumes}, Szarek and Voiculescu proved the free entropy power inequality: for bounded self-adjoint random variables $X_1,\ldots,X_n$ that are free,
\begin{equation}\label{eq:free entropy power}
\exp(2\Sigma(X_1+\cdots+X_n))
\geq \sum_{i=1}^n \exp(2\Sigma(X_i)).
\end{equation}
Shlyakhtenko \cite{Shlyakhtenko2007} later proved a free analogue of \eqref{eq:leave one out entropy power}, from which the monotonicity of free entropy along the free central limit theorem follows. Entropy convergence for normalized free sums was subsequently established by Wang \cite{Wang2010}, while quantitative expansions were obtained by Chistyakov and Götze \cite{ChistyakovGotze2013free} with the same rate $O(n^{-1})$ as the classical case.

Boolean independence is another universal notion of noncommutative independence \cite{bozejko1986positive,schurmann1995direct,Speicher1996}. The Boolean additive convolution was introduced in the work of Speicher and Woroudi \cite{SpeicherWoroudi1997}. In their notation, the self-energy transform plays the role of the $R$-transform in Boolean probability and they also proved a Boolean central limit theorem with the Rademacher law $\mathbf{b}=\frac12(\delta_{-1}+\delta_1)$
being the limiting distribution.

A notion of Boolean entropy was recently introduced in \cite{pan2026boolean}. For a bounded self-adjoint random variable $X$ with distribution $\mu$, it is the logarithmic functional
\[
 \Gamma(X)=\Gamma(\mu)=\int_{\R}\log|x|\,d\mu(x)\in [-\infty,\infty).
\]
Its form is motivated by large-deviation principles for random matrix models in which Boolean independence appears asymptotically. The functional $\Gamma$ has several properties analogous to entropy in the classical and free settings. For a bounded self-adjoint random variable $X$ with a symmetric law, C\'ebron and the author introduced the corresponding Boolean Fisher information \cite{CebronPan2026}. In the present one-dimensional setting it is given by
\[
 \Psi(X)=\Psi(\mu)=\iint_{\R^2}
 \frac{\log x^2-\log y^2}{x^2-y^2}
 \,d\mu(x)\,d\mu(y)\in (0,\infty],
\]
with the diagonal value understood by continuous extension. They proved a Boolean de Bruijn identity and several functional inequalities. In particular, they proved a Boolean entropy power inequality: for a sequence of bounded, self-adjoint and symmetric random variables $X_1,\ldots,X_n$ that are Boolean independent, one has
\begin{equation}\label{eq:boolean entropy power}
    \exp(2\Gamma(X_1+\cdots+X_n))\ge \sum_{i=1}^n \exp(2\Gamma(X_i)).
\end{equation}
Although in \cite{pan2026boolean} the monotonicity of Boolean entropy in the Boolean central limit theorem was proved, the approach is different from the proofs in the classical and free cases. In particular, a Boolean analogue of \eqref{eq:leave one out entropy power} is still unknown. Moreover, the convergence of Boolean entropy in the Boolean central limit theorem is still open, and, if it is true, whether one can say something about the rate of convergence? The main purpose of this article is to answer these questions.

In \cite{CebronPan2026}, the Boolean entropy power inequality \eqref{eq:boolean entropy power} was obtained by first proving the Stam inequality
\begin{equation}\label{eq:stam ineq}
\frac{1}{\Psi(X_1+\cdots+X_n)}
\geq
\sum_{i=1}^n \frac{1}{\Psi(X_i)},
\end{equation}
and then combining it with the Boolean de Bruijn identity. This approach is analogous to the classical and free settings, as Fisher information is more conveniently studied through the score function or its noncommutative analogue, the conjugate variable. It is therefore natural to ask whether the Boolean Stam inequality admits a corresponding extension to arbitrary families of partial sums:
$$
\frac{1}{\Psi(X_1+\cdots+X_n)}\ge \frac{1}{r}\sum_{A\in\mathcal{C}} \frac{1}{\Psi(\sum_{i\in A }X_i)},
$$
which could imply a Boolean entropy power inequality as in \eqref{eq:fractional entropy power} using the de Bruijn identity.

However, such an extension is currently beyond reach. Instead, we establish the following alternative extension of Stam's inequality. We emphasize that the resulting inequality is neither weaker nor stronger than the desired form stated above.

\begin{theorem}\label{thm:fractional-stam}
Let $X_1,\ldots,X_n$ be Boolean independent self-adjoint random variables with symmetric laws and let $[n]=\{1,\ldots,n\}.$ For any $A\subseteq[n]$, set $S_A=\sum_{i\in A}X_i.$ Let $\mathcal C$ be a collection of nonempty subsets of $[n]$, and let $(\beta_A)_{A\in\mathcal C}$ be a fractional partition, namely
$$
\beta_A\geq 0,
\qquad
\sum_{A\in\mathcal C:\,i\in A}\beta_A=1
\quad\text{for every }i=1,\ldots,n.
$$
Then
\begin{equation}\label{eq:fractional stam}
    \frac{1}{\Psi(S_{[n]})}\geq\sum_{A\in\mathcal C}\frac{\beta_A}{\Psi(S_A)}.
\end{equation}

\end{theorem}
We use Theorem~\ref{thm:fractional-stam} in Section~\ref{sec:fractional} to derive some entropy and information inequalities. In particular, when $\mathcal{C}$ is taken to be the subsets of $n-1$ elements and $\beta_A=1/(n-1)$, we obtain an entropy power inequality of the form \eqref{eq:leave one out entropy power} through the Boolean de Bruijn identity. As a consequence, this provides an information-theoretic route to the monotonicity of Boolean entropy along the Boolean central limit theorem in the bounded symmetric case, and the general case was previously established in \cite{pan2026boolean} by a different analytical argument. 

In the i.i.d. setting, let \(\mu_n\) denote the law of the normalized sum
$$
\frac{X_1+\cdots+X_n}{\sqrt n},
$$
where \(X_1,X_2,\ldots\) are identically distributed and Boolean independent self-adjoint random variables with distribution $\mu$. If $\mu$ is centered with variance $1$, then \(\mu_n\) converges weakly to the Rademacher distribution \(\mathbf b\). Moreover, with respect to the weak topology, the sequence \((\mu_n)_{n\geq1}\) admits a continuous extension \((\mu_t)_{t\geq1}\) (see Section \ref{sec:preliminaries} for the exact definition), which preserves the first two moments and has the same limiting law. The following theorem gives a necessary and sufficient condition for entropic convergence in the Boolean central limit theorem.
\begin{theorem}\label{thm:entropy-convergence}
Let $\mu$ and $\mu_t$ be as defined above. If $\Gamma(\mu)>-\infty$, then
\[
 \Gamma(\mu_t)\longrightarrow\Gamma(\mathbf{b})=0.
\]
If $\Gamma(\mu)=-\infty$, then $\Gamma(\mu_t)=-\infty$ for every $t\ge1$.
\end{theorem}

Also, we establish some sharp quantitative asymptotics for the rate of convergence under the finite fourth moment. The next theorem gives the leading order of the relative Boolean entropy. Its coefficient is naturally expressed in terms of the moments of the initial law. However, unlike in the classical and free settings, the rate we get is $O(n^{-1}\log (n))$ instead of $O(n^{-1})$. Example \ref{ex:four-point}, given in Section \ref{sec:entropy-clt}, attains this order.

\begin{theorem}\label{thm:quantitative-entropy}
If $m_4(\mu)<\infty$ and $\Gamma(\mu)>-\infty$,
then
\begin{equation}\label{eq:asymp entropic clt}
    \Gamma(\mu_n\mid \mathbf{b})
 =\frac{m_4(\mu)-1-m_3(\mu)^2}{2}\frac{\log n}{n}
 +o\!\left(\frac{\log n}{n}\right).
\end{equation}
\end{theorem}
Finally, we study the convergence rate of the Fisher information along the Boolean central limit theorem. We emphasize that the Boolean Fisher information is defined only for symmetric probability measures, which is why the third moment term disappears in the coefficient.
\begin{theorem}\label{thm:fisher-asymptotic}
If $\mu$ is symmetric, $m_4(\mu)<\infty$ and $\Psi(\mu)<\infty$, then
\begin{equation}\label{eq:asymp information clt}
    \Psi(\mu_n\mid \mathbf{b})
 =2\bigl(m_4(\mu)-1\bigr)\frac{\log n}{n}
 +o\!\left(\frac{\log n}{n}\right).
\end{equation}
\end{theorem}

\noindent\textbf{Organization of the paper.} The paper is organized as follows. Section~\ref{sec:preliminaries} recalls basic facts about Boolean independence and additive convolution, together with the relevant analytical tools. We also introduce Boolean entropy and Fisher information, and collect several useful representations and lemmas that will be used throughout the paper. Section~\ref{sec:fractional} proves Theorem~\ref{thm:fractional-stam} and derives its applications to entropy and information inequalities. Section~\ref{sec:entropy-clt} is devoted to the proofs of Theorems~\ref{thm:entropy-convergence} and~\ref{thm:quantitative-entropy}, as well as an explicit example demonstrating the sharpness of the quantitative estimate. In the same section, we also prove Theorem~\ref{thm:fisher-asymptotic}.

\section{Preliminaries}\label{sec:preliminaries}

\subsection{Boolean independence and additive convolution}
Consider a von Neumann algebra $\mathcal{M}$ acting on a Hilbert space $\mathcal{H}$, equipped with a vector state
$\varphi(\,\cdot\,)=\langle\,\cdot\,\xi,\xi\rangle,$ where $\xi\in\mathcal{H}$ is a unit vector. In such a noncommutative probability space, random variables are elements in $\mathcal{M}$, and the distribution of a random variable $X\in\mathcal{M}$ is determined by its $*$-moments
\[
\big\{\varphi(P(X,X^*))\big\}_{P\in\mathbb{C}\langle x,y\rangle},
\]
where $\mathbb{C}\langle x,y\rangle$ denotes the algebra of noncommutative polynomials in two indeterminates. For a self-adjoint element $X=X^*$, these moments determine a unique compactly supported probability measure on $\mathbb{R}$, which we denote by $\mathcal{L}(X)$.

We recall the notion of Boolean independence from \cite{bozejko1986positive}. A family of $*$-subalgebras $(\mathscr{A}_i)_{i\in I}\subset\mathcal{M}$ is said to be \emph{Boolean independent} with respect to $\varphi$ if, for every $n\geq1$, every choice of indices $i_1,\ldots,i_n\in I$ satisfying $i_k\neq i_{k+1},\  1\leq k<n,$
and every $X_k\in\mathscr{A}_{i_k}$, one has
\[
\varphi(X_1X_2\cdots X_n)=\varphi(X_1)\varphi(X_2)\cdots\varphi(X_n).
\]
The reader should keep in mind that $\mathrm{Id}$ is not generally Boolean independent from other operators. Indeed, if $X$ is Boolean independent with $\mathrm{Id}$, this implies a trivial factorization: $\varphi(X^n)=\varphi(X)^n$. Therefore, the algebras $\mathscr{A}_i$ are not required to be unital.

We say that $X,Y\in\mathcal{M}$ are \emph{Boolean independent} if the $*$-algebras
\[
\operatorname{Alg}(X)
=
\{P(X,X^*):P\in\mathbb{C}\langle x,y\rangle,\ P(0)=0\}
\quad\text{and}\quad
\operatorname{Alg}(Y)
\]
are Boolean independent.

In particular, if $X,Y\in\mathcal{M}$ are self-adjoint and Boolean independent, with laws $\mu$ and $\nu$, respectively, then the \emph{Boolean additive convolution} $\mu\uplus\nu$ is defined as the law of $X+Y$:
\[
\mathcal{L}(X+Y)=\mu\uplus\nu.
\]
By the Boolean moments-cumulants formula, this convolution admits a convenient analytic description in terms of Cauchy transforms \cite{SpeicherWoroudi1997}. The same analytic definition extends to probability measures with unbounded support.

For completeness, we also recall the corresponding framework for unbounded operators. Let $X$ be an unbounded self-adjoint operator on $\mathcal{H}$. By the spectral theorem, there exists a unique probability measure $\mu\in\mathcal{P}(\mathbb{R})$ such that
\[
\varphi(f(X))=\int_{\mathbb{R}}f(x)\,d\mu(x),
\qquad f\in C_b(\mathbb{R}),
\]
and we call $\mu$ the law of $X$. In \cite{franz2009monotone}, Franz considered the following $C^*$-algebra generated by $X$:
\[
\mathscr{A}_X := \{f(X):f\in C_b(\mathbb{R}), f(0)=0\}.
\]
The vanishing condition is imposed to work with the nonunital algebra generated by $X$, as Boolean independence is naturally formulated for nonunital algebras. Two unbounded self-adjoint operators $X,Y$ are said to be \emph{Boolean independent} if the $C^*$-algebras $\mathscr{A}_X$ and $\mathscr{A}_Y$ are Boolean independent with respect to $\varphi$. With this definition, Franz further proved that if $X$ and $Y$ are Boolean independent self-adjoint operators, then $X+Y$ is essentially self-adjoint and the law of the closure of $X+Y$ is given by $\mathcal{L}(X)\uplus\mathcal{L}(Y)$.

Next, we shall recall the analytical descriptions of Boolean convolution. Let $\Pp(\R)$ denote the space of Borel probability measures on $\R$, and let $\Psym(\R)$ denote the subspace of symmetric probability measures. For $\mu\in\Pp(\R)$, we write

$$
m_k(\mu):=\int_{\R}x^k\,d\mu(x),
$$
whenever the integral exists. We use the notation
$$
\Pzeroone(\R):=\{\mu\in\Pp(\R):m_1(\mu)=0,\ m_2(\mu)=1\}.
$$
For any $c\ne 0$, let $D_c\mu$ denote the pushforward of $\mu$ under the dilation $x\mapsto cx$.

For $\mu\in\Pp(\R)$, the associated Cauchy transform, reciprocal Cauchy transform, and Boolean self-energy transform are defined by
$$
G_\mu(z)=\int_{\R}\frac{d\mu(x)}{z-x},
\qquad
F_\mu(z)=\frac1{G_\mu(z)},
\qquad
K_\mu(z)=z-F_\mu(z),
\qquad z\in\Cplus.
$$
Boolean additive convolution is characterized by the linearization formula
$$
K_{\mu\uplus\nu}=K_\mu+K_\nu.
$$
For $t\ge1$, the Boolean convolution power $\mu^{\uplus t}$ is defined by
$$
K_{\mu^{\uplus t}}=tK_\mu.
$$
The transforms behave naturally under dilation:
$$
G_{D_c\mu}(z)=c^{-1}G_\mu(c^{-1}z),
\qquad
K_{D_c\mu}(z)=cK_\mu(c^{-1}z).
$$

For $\mu\in\Pp(\R)$, let $\mu^{(2)}$ denote the pushforward of $\mu$ under the map $x\mapsto x^2$. If $X_1,\ldots,X_m$ are Boolean independent self-adjoint random variables with symmetric laws, then
\begin{equation}\label{eq:square-convolution}
\mathcal L\left(\left(\sum_{i=1}^m X_i\right)^2\right)
=\biguplus_{i=1}^m \mathcal L(X_i^2).
\end{equation}
Equivalently, if $\mu$ is symmetric and $\rho=\mu^{(2)}$, then

$$
G_\mu(z)=zG_\rho(z^2),
\qquad
F_\mu(z)=\frac{F_\rho(z^2)}{z},
\qquad
K_\mu(z)=\frac{K_\rho(z^2)}{z}.
$$
In the setting of the Boolean central limit theorem, for $\mu\in\Pzeroone(\R)$, we define the sequence of laws of the normalized sums in the Boolean central limit theorem”:
$$
\mu_n
:=
\underbrace{D_{\frac{1}{\sqrt n}}\mu\uplus\cdots\uplus
D_{\frac{1}{\sqrt n}}\mu}_{n\ \mathrm{times}},
\qquad n\in\mathbb N.
$$
Their Cauchy transforms are given by
\begin{equation*}
G_{\mu_n}(z)
=\frac{1}{z-\sqrt{n}K_\mu(\sqrt{n}z)},\qquad n\in \mathbb N.
\end{equation*}
More generally, the parameter $n$ can be replaced by any real number $t\ge1$. This gives a continuous family of probability measures $(\mu_t)_{t\ge1}\subset\Pzeroone(\R)$, which extends the sequence $(\mu_n)_{n\ge1}$. The family satisfies a multiplicative semigroup property
\begin{equation}\label{eq:semigroup}
    (\mu_t)_s=\mu_{ts},\qquad \forall s,t\ge1,
\end{equation}
see \cite[Lemma 4.2]{pan2026boolean} for more details.

\subsection{Nevanlinna representation}\label{sec:nevanlinna}

We recall a standard representation: if $\varphi:\mathbb{C}_+\to \mathbb{C}_+$ is analytic, then there exists a unique finite positive Borel measure $\sigma$ on $\mathbb{R}$ and real numbers $\alpha,\beta$ with $\beta \ge 0$ such that
\begin{equation*}
    \varphi(z)=\alpha +\beta z +\int\frac{1+tz}{t-z}\,d\sigma(t).
\end{equation*}
If, in addition, $\lim_{y\to \infty}y|\varphi(iy)|=c<\infty,$ then the above representation takes the simple form:
$$
\varphi(z) = \int \frac{1}{t-z}\,d\tilde{\sigma}(t) = -G_{\tilde{\sigma}}(z),
$$
where $d\tilde{\sigma}(t)=(1+t^2)d\sigma(t)$ with total mass $c$ (see e.g. \cite[Theorem 10]{mingo2017free} for more details). Applying this to the self-energy transform $K_\mu$, Speicher and Woroudi \cite{SpeicherWoroudi1997} showed that if $\mu\in \mathcal{P}_{0,1}(\mathbb{R})$, then there exists a probability measure $\nu$ on $\mathbb{R}$ such that
\begin{equation}\label{eq:Nevanlinna repre P_0,1}
    K_\mu(z)=G_\nu(z).
\end{equation}
Equivalently, $F_\mu(z)=z-G_\nu(z).$ The moments of $\nu$ are related to the moments of $\mu$ through the expansion of Cauchy transforms at infinity. If $\mu$ has finite fourth moment, then
$$
G_\mu(z)
=
\frac{1}{z}
+\frac{m_1(\mu)}{z^2}
+\frac{m_2(\mu)}{z^3}
+\frac{m_3(\mu)}{z^4}
+\frac{m_4(\mu)}{z^5}
+o(|z|^{-5}),
$$
and hence, when $m_1(\mu)=0$ and $m_2(\mu)=1$,
$$
F_\mu(z)
=
z-\frac{1}{z}
-\frac{m_3(\mu)}{z^2}
-\frac{m_4(\mu)-1}{z^3}
+o(|z|^{-3}).
$$
Therefore, by definition,
$$
G_\nu(z)
=
\frac{1}{z}
+\frac{m_1(\nu)}{z^2}
+\frac{m_2(\nu)}{z^3}
+o(|z|^{-3}).
$$
Comparing coefficients yields
$$
m_1(\nu)=m_3(\mu),
\qquad
m_2(\nu)=m_4(\mu)-1.
$$
In particular,
$$
m_4(\mu)-1-m_3(\mu)^2=\operatorname{Var}(\nu)\geq0.
$$
This variance will be the coefficient governing the leading logarithmic correction in the quantitative Boolean central limit theorem, see \eqref{eq:asymp entropic clt} and \eqref{eq:asymp information clt}.

\subsection{Boolean entropy and Fisher information}
We recall the definitions and some basics about Boolean entropy and Fisher information from \cite{pan2026boolean,CebronPan2026}. For any $\mu\in\Pp(\R)$ with finite second moment, the \emph{Boolean entropy} is given by
\[
 \Gamma(\mu)=\int_{\R}\log|x|\,d\mu(x)\in[-\infty,\infty).
\]
The positive part of the logarithm is integrable under the second-moment assumption. In the space $\mathcal{P}_{0,1}(\mathbb{R})$, the Boolean entropy is maximized by the Rademacher distribution $\mathbf{b}$. If $\mu$ is symmetric, the \emph{Boolean Fisher information} is
\[
 \Psi(\mu)=\iint_{\R^2}
 \frac{\log x^2-\log y^2}{x^2-y^2}
 \,d\mu(x)\,d\mu(y),
\]
where the diagonal value is the continuous extension $x^{-2}$. In terms of the square pushforward,
\[
 \Psi(\mu)=\Psi^{(2)}(\mu^{(2)}),
 \qquad
 \Psi^{(2)}(\rho)=\int_0^\infty\int_0^\infty
 \frac{\log x-\log y}{x-y}
 \,d\rho(x)\,d\rho(y),\quad \rho\in \mathcal{P}(\mathbb{R}_+).
\]
The scaling relations are
\[
 \Gamma(\D_c\mu)=\Gamma(\mu)+\log |c|,
 \qquad
 \Psi(\D_c\mu)=c^{-2}\Psi(\mu),\quad c\ne 0.
\]

For a self-adjoint operator $X$ with distribution $\mu$, we write $\Gamma(X)=\Gamma(\mu),$ and $\Psi(X)=\Psi(\mu)$ if $\mu$ is symmetric. Let $B$ be a Rademacher-distributed variable ($\mathcal{L}(B)=\mathbf{b}$), Boolean independent from $X$. If $X$ is bounded such that $\Gamma(X)>-\infty$ and $\Psi(X)<\infty$, then the following \emph{Boolean de Bruijn identity} holds,
\begin{equation}\label{eq:debruijn}
    \Gamma(X)=\frac12\int_0^\infty
 \left(\frac1{1+s}-\Psi(X+\sqrt{s}\,B)\right)ds.
\end{equation}

For any $\mu\in\mathcal{P}(\mathbb{R})$ with finite variance, the \emph{Boolean entropy and Fisher information relative to the Rademacher distribution} are defined as
\begin{align}
    \Gamma(X\mid B)=\Gamma(\mu\mid \mathbf{b}) &= \Gamma(\mathbf{b})-\frac{1}{2}m_2(\mathbf{b})-\Gamma(\mu)+\frac{1}{2}m_2(\mu)\ge 0;\label{eq:relative entropy}\\
    \Psi(X\mid B)=\Psi(\mu\mid \mathbf{b}) &= \Psi(\mu)+m_2(\mu)-\Psi(\mathbf{b})-m_2(\mathbf{b})\ge 0\label{eq:relative fisher}.
\end{align}
In particular, when $\mu$ has unit second moment, note that $\Gamma(\mathbf{b})=0,\,\Psi(\mathbf{b})=1,$ so that we have $\Gamma(\mu\mid \mathbf{b})=-\Gamma(\mu),\,\Psi(\mu\mid\mathbf{b})=\Psi(\mu)-1$.

\subsection{A negative-axis representation}\label{sec:negative-axis}
For any $\rho\in\Pp(\mathbb{R}_+)$ and $s>0$, set
\[
 h_\rho(s)=-G_\rho(-s)=\int_0^\infty\frac{d\rho(x)}{x+s},
 \qquad
 r_\rho(s)=K_\rho(-s)=\frac1{h_\rho(s)}-s.
\]
\begin{lemma}\label{lem:stieltjes}
For every $\rho\in\Pp(\mathbb{R}_+)$,
\begin{equation}\label{eq:stieltjes}
    \Psi^{(2)}(\rho)
 =\int_0^\infty h_\rho(s)^2\,ds
 =\int_0^\infty\frac{ds}{(s+r_\rho(s))^2}.
\end{equation}
Moreover, $r_\rho$ is nonnegative and nondecreasing on $(0,\infty)$.
\end{lemma}

\begin{proof}
For $x,y>0$, the elementary identity
\[
 \frac{\log x-\log y}{x-y}
 =\int_0^\infty\frac{ds}{(x+s)(y+s)}
\]
allows us, by Fubini's theorem, to obtain the first equality in \eqref{eq:stieltjes} and hence the representation in terms of $r_\rho$. Since $h_\rho(s)\le s^{-1}$, we have $r_\rho(s)\ge0$. Finally, differentiation gives
\[
 r_\rho'(s)
 =\frac{\displaystyle\int (x+s)^{-2}\,d\rho(x)}{\displaystyle\left(\int (x+s)^{-1}\,d\rho(x)\right)^2}-1\ge0
\]
by Cauchy--Schwarz.
\end{proof}

\begin{lemma}\label{lem:boolean-power}
For any $\rho\in\Pp(\mathbb{R}_+)$ and $c\ge1$, we have
\begin{equation}
    \Psi^{(2)}(\rho^{\uplus c})\le\frac1c\Psi^{(2)}(\rho).
\end{equation}
\end{lemma}

\begin{proof}
Since $K_{\rho^{\uplus c}}=cK_\rho$, we have $r_{\rho^{\uplus c}}(s)=cr_\rho(s)$. Therefore, by Lemma \ref{lem:stieltjes} and the change of variables $s=ca$,
\[
 \Psi^{(2)}(\rho^{\uplus c})
 =\int_0^\infty\frac{ds}{(s+cr_\rho(s))^2}
 =\frac1c\int_0^\infty\frac{da}{(a+r_\rho(ca))^2}.
\]
The monotonicity of $r_\rho$ gives $r_\rho(ca)\ge r_\rho(a)$, proving the claim.
\end{proof}

\begin{lemma}\label{lem:representation relative entropy}
    For any $\mu\in\mathcal{P}(\mathbb{R})$ with $m_2(\mu)=1$ and $\Gamma(\mu)>-\infty$, we have
    \begin{equation}\label{eq:relative-entropy-integral}
        \Gamma(\mu\mid \mathbf{b}) = \frac{1}{2}\int_0^\infty \left(h_{\mu^{(2)}}(s)-\frac{1}{1+s}\right)\,ds
    \end{equation}
\end{lemma}
\begin{proof}
Recall that if $m_2(\mu)=1$, one has
$$
\Gamma(\mu\mid\mathbf b) =-\Gamma(\mu),
\qquad
\Gamma(\mu)=\int_{\mathbb R}\log|x|\,d\mu(x),
$$
For every $x\neq0$, note that the following identity holds
$$
\log|x|
=
\frac12\int_0^\infty
\left(\frac{1}{1+s}-\frac{1}{x^2+s}\right)\,ds,
$$
which follows from
$$
\int_0^\infty
\left(\frac{1}{s+a}-\frac{1}{s+b}\right)\,ds
=\log\frac{b}{a},
\qquad a,b>0.
$$
Therefore,
$$
\Gamma(\mu)
=
\frac12\int_{-\infty}^{+\infty}\int_0^\infty
\left(\frac{1}{1+s}-\frac{1}{x^2+s}\right)
\,ds\,d\mu(x).
$$
Since $\Gamma(\mu)>-\infty$, the interchange of the two
integrals is justified by splitting the integral into the regions
$|x|\leq1$ and $|x|>1$. Indeed,
$$
\int_0^\infty
\left|
\frac{1}{1+s}-\frac{1}{x^2+s}
\right|\,ds
=2|\log|x||,
$$
and $\log^+|x|$ is integrable whenever $\mu$ has finite second
moment, while $\Gamma(\mu)>-\infty$ guarantees the integrability
of $\log^-|x|$.
Hence
$$
\Gamma(\mu)
=
\frac12\int_0^\infty
\left(
\frac{1}{1+s}
-
\int_{\mathbb R}\frac{1}{x^2+s}\,d\mu(x)
\right)\,ds,
$$
which yields the desired formula.
\end{proof}
\begin{remark}
    The representation in Lemma \ref{lem:representation relative entropy} admits a useful interpretation in terms of the nonmicrostates Boolean Fisher information. Recall that for a probability measure $\mu$, the nonmicrostates Fisher information of $\mu$ is defined as
    $$
    \Psi^*(\mu)=\int_{\mathbb{R}} \frac{1}{x^2}\,d\mu(x).
    $$
    Now let $(\tilde{\mu}_t)_{t\geq0}$ be the flow defined by
    $$
    \tilde{\mu}_t=(T_t)_\#\mu,
    \quad T_t(x)=\operatorname{sgn}(x)\sqrt{x^2+t}.
    $$
    That is, $\tilde{\mu}_t^{(2)}=(x\mapsto x+t)_\#\mu^{(2)},$ and hence
    $$
    \Psi^*(\tilde{\mu}_t)=\int_{\mathbb R}\frac{1}{x^2}\,d\tilde{\mu}_t(x)=\int_0^\infty \frac{1}{x+t}\,d\mu^{(2)}(x)=h_{\mu^{(2)}}(t).
    $$
    Moreover,
    $$
    \frac{d}{dt}\Gamma(\tilde{\mu}_t)=\frac12\Psi^*(\tilde{\mu}_t).
    $$
    Indeed, $T_t'(x)=1/(2T_t(x))$, and therefore
    $$
    \frac{d}{dt}\log|T_t(x)|=\frac{1}{2T_t(x)^2}.
    $$
    This should be compared with the nonmicrostates de Bruijn identity along the Boolean heat flow established in~\cite{CebronPan2026}.
\end{remark}

\section{A fractional Stam inequality}\label{sec:fractional}

Let $n\ge1$, and write $[n]=\{1,\ldots,n\}$. A family $(\beta_A)_{A\in\mathcal C}$ indexed by a collection of nonempty subsets of $[n]$ is a fractional partition if
\[
 \beta_A\ge0,
 \qquad
 \sum_{A\in\mathcal C:\,i\in A}\beta_A=1
 \quad\text{for every }i=1,\ldots,n.
\]

We now prove Theorem~\ref{thm:fractional-stam}.

\begin{proof}
First assume that $\Psi(X_i)<\infty$ for all $i=1,\ldots,n.$ By the Stam inequality \eqref{eq:stam ineq}, we have $\Psi(S_A)<\infty$ for all $A\subset [n]$. Set $\rho_i=\law(X_i^2),$ $\rho_A=\law(S_A^2)$. By \eqref{eq:square-convolution},
\begin{equation}\label{eq:r-additive}
    K_{\rho_A}=\sum_{i\in A}K_{\rho_i},
 \qquad
 r_{\rho_A}(s)=\sum_{i\in A}r_{\rho_i}(s).
\end{equation}
Thus Lemma \ref{lem:stieltjes} gives
\[
 \Psi(S_A)=\int_0^\infty
 \frac{ds}{\left(s+\sum_{i\in A}r_{\rho_i}(s)\right)^2}.
\]
Since each $r_{\rho_i}$ is nonnegative,
\begin{equation}\label{eq:monotone-psi}
    A\subseteq B\quad\Longrightarrow\quad \Psi(S_B)\le\Psi(S_A).
\end{equation}
Set
\[
 a_A=\frac1{\Psi(S_A)}>0,
 \qquad
 T=\sum_{A\in\mathcal C}\beta_Aa_A.
\]
If $a_A\ge T$ for some $A$, then \eqref{eq:monotone-psi} immediately gives
\[
 \frac1{\Psi(S_{[n]})}\ge a_A\ge T.
\]
We may therefore assume $a_A<T$ for every $A$. Define
\[
 c_A=\frac{T}{a_A}>1,
 \qquad
 w_A=\frac{\beta_Aa_A}{T}.
\]
Then $\sum_Aw_A=1$ and $w_Ac_A=\beta_A$.
Let $\tilde\rho_A=\rho_A^{\uplus c_A}$. By the definition of fractional partition and \eqref{eq:r-additive},
\[
 \sum_Aw_AK_{\tilde\rho_A}
 =\sum_A\beta_AK_{\rho_A}
 =\sum_{i=1}^m\left(\sum_{A\ni i}\beta_A\right)K_{\rho_i}
 =K_{\rho_{[n]}}.
\]
Since $K_\rho(z)=z-F_\rho(z)$ and $\sum_A w_A=1$, we obtain
\[
F_{\rho_{[n]}}(z)=\sum_A w_A F_{\tilde\rho_A}(z).
\]
Evaluating at $z=-s$ for $s>0$, we have
$F_{\rho_{[n]}}(-s)<0$ and $F_{\tilde\rho_A}(-s)<0$ for every $A$. Since $x\mapsto x^{-2}$ is convex on $(-\infty,0)$, Jensen's inequality therefore gives
\[
 G_{\rho_{[n]}}(-s)^2
 \le\sum_Aw_AG_{\tilde\rho_A}(-s)^2.
\]
Integrating with respect to $s$ and applying Lemma \ref{lem:stieltjes} and Lemma \ref{lem:boolean-power},
\[
 \Psi(S_{[n]})
 \le\sum_Aw_A\Psi^{(2)}(\tilde\rho_A)
 \le\sum_A\frac{w_A}{c_A}\Psi^{(2)}(\rho_A)
 =\frac1T.
\]
Taking reciprocals yields the desired inequality.

If \(\Psi(X_i)=\infty\) for some \(i\in[n]\), we regularize the corresponding law by setting \(\rho_i(\varepsilon):=\rho_i\uplus\delta_\varepsilon\). Then \(r_{\rho_i(\varepsilon)}(s)=r_{\rho_i}(s+\varepsilon)\). Applying the previous argument to the perturbed laws and then by Lemma \ref{lem:stieltjes}, letting \(\varepsilon\downarrow0\), the monotone convergence theorem completes the proof.
\end{proof}

Take $A_j=[n]\setminus\{j\}$. Each index belongs to exactly $n-1$ of these sets, so the weights $\beta_{A_j}=1/(n-1)$ form a fractional partition. Applying Theorem \ref{thm:fractional-stam} gives the following result.
\begin{corollary}\label{cor:leave-one-out}
Let $X_1,\ldots,X_{n}$ be Boolean independent self-adjoint random variables with symmetric laws. Then
\[
 \frac{1}{\Psi\!\left(\sum_{i=1}^{n}X_i\right)}
 \ge \frac{1}{n-1}\sum_{j=1}^{n}
 \frac{1}{\Psi\!\left(\sum_{i\ne j}X_i\right)}.
\]
\end{corollary}

We also deduce the improved Stam inequality for weighted sums.
\begin{corollary}\label{cor:weighted-stam}
Under the assumptions of Corollary \ref{cor:leave-one-out}, let $a=(a_1,\ldots,a_{n})\in\R^{n}$ satisfy $\sum_i a_i^2=1$, and set $S_a=\sum_i a_iX_i$. For $|a_j|<1$, define
\[
 S_a^{(j)}=\frac1{\sqrt{1-a_j^2}}\sum_{i\ne j}a_iX_i.
\]
Then
\[
 \frac1{\Psi(S_a)}
 \ge\frac{1}{n-1}\sum_{j=1}^{n}\frac{1-a_j^2}{\Psi(S_a^{(j)})},
\]
where terms with $a_j^2=1$ are omitted, and
\begin{equation}\label{eq:weighted stam}
    \Psi(S_a)\le\frac{1}{n-1}\sum_{j=1}^{n}(1-a_j^2)\Psi(S_a^{(j)}).
\end{equation}
\end{corollary}

\begin{proof}
Apply Corollary \ref{cor:leave-one-out} to $a_1X_1,\ldots,a_{n}X_{n}$. By scaling,
\[
 \frac1{\Psi\left(\sum_{i\ne j}a_iX_i\right)}
 =\frac{1-a_j^2}{\Psi(S_a^{(j)})},
\]
which gives the reciprocal inequality. Setting $p_j=(1-a_j^2)/(n-1)$, we have $p_j\ge0$ and $\sum_jp_j=1$. The weighted harmonic-arithmetic mean inequality then yields the second bound.
\end{proof}

\begin{corollary}\label{cor:weighted-entropy}
Under the assumptions of Corollary \ref{cor:weighted-stam}, assume further that $X_1,\ldots,X_n$ are bounded, then we have
\begin{equation}\label{eq:weighted entropy ineq}
    \Gamma(S_a)
 \ge\frac{1}{n-1}\sum_{j=1}^{n}(1-a_j^2)\Gamma(S_a^{(j)}).
\end{equation}
\end{corollary}

\begin{proof}
Let $B_1,\ldots,B_{n}$ be Boolean independent Rademacher variables independent of the $X_i$. The square-convolution identity shows that both
\[
 \sum_{i=1}^{n}a_iB_i
 \quad\text{and}\quad
 \frac1{\sqrt{1-a_j^2}}\sum_{i\ne j}a_iB_i
\]
have the Rademacher law. Applying the weighted Stam inequality \eqref{eq:weighted stam} to $X_i+\sqrt s B_i$ gives, for $s>0$,
\begin{equation}\label{eq:weighted-fisher-debruijn}
    \Psi(S_a+\sqrt s B)
 \le\frac{1}{n-1}\sum_{j=1}^{n}(1-a_j^2)
 \Psi(S_a^{(j)}+\sqrt s B_j').
\end{equation}
Subtracting from $(1+s)^{-1}$, integrating in $s$, and using the Boolean de Bruijn identity \eqref{eq:debruijn} yields the stated entropy inequality.
\end{proof}
The boundedness assumption is required because of the de Bruijn identity \eqref{eq:debruijn}. As a corollary, we obtain the Boolean analogue of \eqref{eq:leave one out entropy power}.
\begin{corollary}\label{cor:epine-many}
Under the assumptions of Corollary \ref{cor:weighted-entropy}, set $S_n=\sum_{i=1}^{n}X_i$ and $S_n^{(j)}=\sum_{i\ne j}X_i$. Then
\begin{equation}
    \exp\bigl(2\Gamma(S_n)\bigr)
 \ge\frac{1}{n-1}\sum_{j=1}^{n}
 \exp\bigl(2\Gamma(S_n^{(j)})\bigr).
\end{equation}
\end{corollary}

\begin{proof}
Set
\[
 N_j=e^{2\Gamma(S_n^{(j)})},
 \qquad
 N=\frac{1}{n-1}\sum_{j=1}^{n}N_j.
\]
If $N_j\ge N$ for some $j$, then \eqref{eq:boolean entropy power} gives $e^{2\Gamma(S_n)}\ge N$. Otherwise define $q_j=N_j/N$ and $a_j=\sqrt{1-q_j}$. Since $\sum_jq_j=n-1$, we have $\sum_ja_j^2=1$. Applying Corollary \ref{cor:weighted-entropy} to the rescaled variables $a_j^{-1}X_j$ gives
\[
 \Gamma(S_n)
 \ge\frac1{n-1}\sum_{j=1}^{n}q_j
 \Gamma\bigl(q_j^{-1/2}S_n^{(j)}\bigr).
\]
Using scaling once more,
\[
 \Gamma(S_n)
 \ge\frac1{n-1}\sum_jq_j
 \left(\Gamma(S_n^{(j)})-\frac12\log q_j\right)
 =\frac12\log N.
\]
Exponentiating proves the result.
\end{proof}

\section{The entropic Boolean central limit theorem}\label{sec:entropy-clt}

\subsection{Entropy convergence in the Boolean CLT}
For $\mu\in\Pzeroone(\R)$, recall the Boolean central limit flow $\mu_t=\D_{t^{-1/2}}(\mu^{\uplus t}),\, t\ge1$ is defined via its Cauchy transform
\begin{equation}\label{eq:flow-transform}
G_{\mu_t}(z)
=\frac{1}{z-\sqrt{t}K_\mu(\sqrt{t}z)},\quad K_\mu(z)=z-\frac{1}{G_\mu(z)}.
\end{equation}
We now prove Theorem~\ref{thm:entropy-convergence}.
\begin{proof}
For any $s>0$, set
\[
 h_t(s)=\int_{\R}\frac{d\mu_t(x)}{x^2+s}=-\mathrm{Im} \frac{1}{\sqrt{s}}G_{\mu_t}(i\sqrt{s}).
\]
Fix $s>0$. Differentiating with respect to $t$ and using the semigroup relation \eqref{eq:semigroup} gives
\begin{equation}\label{eq:ht-derivative}
    t\partial_th_t(s)
 =-\bigl(\mathrm{Re}\,G_{\mu_t}(i\sqrt s)\bigr)^2
 -s\left(\int_{\R}\frac{d\mu_t(x)}{(x^2+s)^2}-h_t(s)^2\right).
\end{equation}
The second term on the right-hand side is nonnegative, and we deduce that $t\partial_th_t(s)\le 0$. Moreover,
$$
\operatorname{Re}G_{\mu_t}(i\sqrt s)
=
-\int_{\mathbb R}\frac{x}{x^2+s}\,d\mu_t(x).
$$
Therefore, by Cauchy-Schwarz,
$$
\begin{aligned}
\bigl(\operatorname{Re}G_{\mu_t}(i\sqrt s)\bigr)^2
&\leq
\int_{\mathbb R}\frac{x^2}{(x^2+s)^2}\,d\mu_t(x)\\
&=
\int_{\mathbb R}
\left(
\frac{1}{x^2+s}
-\frac{s}{(x^2+s)^2}
\right)d\mu_t(x)\\
&=
h_t(s)-s\int_{\mathbb R}\frac{d\mu_t(x)}{(x^2+s)^2}.
\end{aligned}
$$
Substituting this inequality into \eqref{eq:ht-derivative} yields
$$
-h_t(s)\le t\partial_t h_t(s)\le0.
$$
The right-hand inequality implies that $t\mapsto h_t(s)$ is decreasing. For the left-hand inequality, we observe that
$
\partial_t\bigl(t h_t(s)\bigr)
=h_t(s)+t\partial_t h_t(s)\ge0.
$
Thus $t\mapsto t h_t(s)$ is increasing, which implies that
$t h_t(s)\ge h_1(s),$ and consequently,
\begin{equation}\label{eq:ht-bounds}
    \frac1t h_1(s)\le h_t(s)\le h_1(s),\quad t\ge1.
\end{equation}

For $\varepsilon>0$, introduce the regularized logarithmic integral
\[
\Gamma_\varepsilon(\nu)=\frac12\int_{\R}\log(x^2+\varepsilon)\,d\nu(x).
\]
Tonelli's theorem gives,
\begin{equation}\label{eq:regularized-entropy}
    \Gamma_\varepsilon(\mu_t)-\Gamma(\mu_t)
 =\frac12\int_0^\varepsilon h_t(s)\,ds.
\end{equation}

Assume first that $\Gamma(\mu)>-\infty$. Taking $t=1$ in \eqref{eq:regularized-entropy} shows that $h_1$ is integrable on every interval $(0,\varepsilon)$. For fixed $s>0$, weak convergence $\mu_t\Rightarrow \mathbf{b}$ gives
\[
 h_t(s)\longrightarrow h_\infty(s)=\frac1{1+s}.
\]
The upper bound in \eqref{eq:ht-bounds} and dominated convergence therefore yield
\begin{equation}\label{eq:regularized-limit}
    \Gamma_\varepsilon(\mu_t)-\Gamma(\mu_t)
 \longrightarrow
 \frac12\int_0^\varepsilon\frac{ds}{1+s}
 =\frac12\log(1+\varepsilon).
\end{equation}
On the other hand, for sufficiently large $M>0$, decompose 
$$
\Gamma_\varepsilon(\mu_t)=\frac{1}{2}\int \log (x^2+\varepsilon)\wedge M\,d\mu_t(x)+\frac{1}{2}\int (\log (x^2+\varepsilon)-\log (x^2+\varepsilon)\wedge M)\,d\mu_t(x).
$$
Since the family $\{\mu_t\}_{t\ge1}$ has uniformly bounded second moments, it is uniformly integrable. Thus, the second term can be made uniformly arbitrarily small. The first term tends to $\Gamma_\varepsilon(\mathbf{b})$ by the weak convergence $\mu_t\Rightarrow \mathbf{b}$. Therefore, letting $M\to \infty$ yields
\[
\Gamma_\varepsilon(\mu_t)\longrightarrow\Gamma_\varepsilon(\mathbf{b})=\frac12\log(1+\varepsilon).
\]
Combining this with \eqref{eq:regularized-limit} proves $\Gamma(\mu_t)\longrightarrow0$.

Now assume $\Gamma(\mu)=-\infty$. Then \eqref{eq:regularized-entropy} at $t=1$ implies
\[
 \int_0^\varepsilon h_1(s)\,ds=\infty.
\]
The lower bound in \eqref{eq:ht-bounds} gives
\[
 \int_0^\varepsilon h_t(s)\,ds
 \ge \frac1t\int_0^\varepsilon h_1(s)\,ds
 =\infty.
\]
Since $\Gamma_\varepsilon(\mu_t)$ is finite, \eqref{eq:regularized-entropy} forces $\Gamma(\mu_t)=-\infty$.
\end{proof}

\subsection{Quantitative entropy asymptotics}
We now prove Theorem~\ref{thm:quantitative-entropy}.
\begin{proof}
By \eqref{eq:Nevanlinna repre P_0,1}, let $\nu$ be such that $F_\mu(z)=z-G_\nu(z),$ and hence
\begin{equation*}
G_n(z):=G_{\mu_n}(z)
=\frac{1}{z-\sqrt n\,G_\nu(\sqrt n z)}.
\end{equation*}
Write $b_3=m_1(\nu)=m_3(\mu),\ b_4=
m_2(\nu)=m_4(\mu)-1$. For any $s>0$, define
\[
R_s:=\int_{\mathbb R}\frac{x}{x^2+s^2}\,d\nu(x),
\qquad
I_s:=\int_{\mathbb R}\frac{1}{x^2+s^2}\,d\nu(x),
\]
so that $G_\nu(is)=-R_s-isI_s.$ By Lemma \ref{lem:representation relative entropy}, after the change of variables
$s\mapsto s^2/n$, \eqref{eq:relative-entropy-integral} becomes
\begin{equation*}
\Gamma(\mu_n\mid\mathbf b)
=
\int_0^\infty
\frac{s}{n}
\left(
f_n(s)-\frac{n}{n+s^2}
\right)\,ds,
\end{equation*}
where
\begin{equation}\label{eq:f_n}
f_n(s)
:=h_n\left(\frac{s^2}{n}\right)=-\operatorname{Im}\frac{1}{s/\sqrt n}
G_n\left(\frac{is}{\sqrt n}\right)=
\frac{I_s+n^{-1}}
{R_s^2+s^2(I_s+n^{-1})^2}.
\end{equation}

\medskip
\noindent\textbf{The compact region.}
Fix $M>0$. By Jensen's inequality,
\[
I_s
\geq \frac{1}{b_4+s^2}
\geq \frac{1}{b_4+M^2},
\qquad 0<s\leq M.
\]
Hence, for all sufficiently large $n$,
$I_s\geq n^{-1}$ on $(0,M]$, and therefore
\[
f_n(s)\leq 2f_\infty(s),
\qquad
f_\infty(s):=
\frac{I_s}{R_s^2+s^2I_s^2}.
\]
On the other hand, $F_\mu(is)
=
is-G_\nu(is)
=
R_s+is(1+I_s),$ so
\begin{equation}\label{eq:ImGmu}
-\operatorname{Im}G_\mu(is)
=
\frac{s(1+I_s)}
{R_s^2+s^2(1+I_s)^2}.
\end{equation}
By the lower bound $I_s\geq (b_4+M^2)^{-1}$, there exists a constant $C_M>1$ such that $1+I_s\leq C_M I_s$ for all $0<s\leq M,$ and consequently $sf_\infty(s)\leq C^2_M\bigl[-\operatorname{Im}G_\mu(is)\bigr].$
Moreover, note that
\[
\begin{aligned}
\int_0^M-\operatorname{Im}G_\mu(is)\,ds
&=
\int_{\mathbb R}
\left(
\int_0^M\frac{s}{x^2+s^2}\,ds
\right)d\mu(x)\\
&=
\frac12
\int_{\mathbb R}
\log\left(1+\frac{M^2}{x^2}\right)d\mu(x)
<\infty.
\end{aligned}
\]
Here the finiteness follows from $\Gamma(\mu)>-\infty$.
Thus
\[
\int_0^M \frac{s}{n} f_n(s)\,ds=O_M(n^{-1}).
\]
Finally, 
\[
\int_0^M
\frac{s}{n}\frac{n}{n+s^2}\,ds
=
\frac12\log\left(1+\frac{M^2}{n}\right)
=
O_M(n^{-1}).
\] 
Therefore
\begin{equation}\label{eq:compact}
\int_0^M
\frac{s}{n}
\left|
f_n(s)-\frac{n}{n+s^2}
\right|\,ds
=
O_M(n^{-1}),
\end{equation}

\medskip
\noindent\textbf{Expansion on the tail.}
We now consider $s\geq M$, and set $A_n(s):=s^{-2}+n^{-1}.$ Since $m_2(\nu)=b_4$,
\[
I_s=s^{-2}-b_4s^{-4}+\varepsilon_s,
\qquad
\varepsilon_s:=
s^{-4}\int_{\mathbb R}\frac{x^4}{x^2+s^2}\,d\nu(x).
\]
Hence
\begin{equation}\label{eq:eps-bound}
0\leq\varepsilon_s\leq b_4s^{-4},
\qquad
s^4\varepsilon_s
=
\int_{\mathbb R}\frac{x^4}{x^2+s^2}\,d\nu(x)
\xrightarrow[s\to\infty]{}0,
\end{equation}
where the convergence follows from the dominated convergence theorem.

Similarly,
\[
R_s=b_3s^{-2}+\eta_s,
\qquad
\eta_s:=
-s^{-2}\int_{\mathbb R}\frac{x^3}{x^2+s^2}\,d\nu(x).
\]
By $|x|/(x^2+s^2)\leq 1/(2s)$, we obtain
\begin{equation}\label{eq:eta-bound}
|\eta_s|
\leq
\frac{b_4}{2}s^{-3}.
\end{equation}

Set
\[
d_s:=I_s-s^{-2}=-b_4s^{-4}+\varepsilon_s.
\]
Then $I_s+n^{-1}=A_n(s)+d_s$, and \eqref{eq:f_n} gives
\[
f_n(s)
=
\frac{s^{-2}(A_n(s)+d_s)}
{s^{-2}R_s^2+(A_n(s)+d_s)^2}.
\]
Let
\[
B_n(s):=s^{-2}R_s^2+(A_n(s)+d_s)^2.
\]
Since $\frac{n}{n+s^2}=\frac{s^{-2}}{A_n(s)},$ we have the exact identity
\begin{equation}\label{eq:exact-difference}
f_n(s)-\frac{s^{-2}}{A_n(s)}
=-\frac{s^{-2}\bigl(A_n(s)d_s+d_s^2+s^{-2}R_s^2\bigr)}
{A_n(s)B_n(s)}.
\end{equation}

Write
\[
R_s^2=b_3^2s^{-4}+\delta_s,
\qquad
\delta_s:=R_s^2-b_3^2s^{-4}.
\]
By \eqref{eq:eta-bound},
\[
|\delta_s|
\leq
2|b_3|s^{-2}|\eta_s|+\eta_s^2
\leq C_\mu s^{-5}.
\]
Substituting
\[
d_s=-b_4s^{-4}+\varepsilon_s,
\qquad
R_s^2=b_3^2s^{-4}+\delta_s
\]
into \eqref{eq:exact-difference}, we obtain
\[
\begin{aligned}
f_n(s)-\frac{s^{-2}}{A_n(s)}
={}&
\frac{b_4s^{-6}}{B_n(s)}
-\frac{b_3^2s^{-8}}{A_n(s)B_n(s)}
-\frac{s^{-2}\varepsilon_s}{B_n(s)}\\
&-\frac{s^{-2}d_s^2}{A_n(s)B_n(s)}
-\frac{s^{-4}\delta_s}{A_n(s)B_n(s)}.
\end{aligned}
\]
Therefore, if we define
\begin{align*}
\xi_n(s):={}&
b_4s^{-6}
\left(\frac{1}{B_n(s)}-\frac{1}{A_n(s)^2}\right)\\
&+b_3^2s^{-8}
\left(\frac{1}{A_n(s)^3}
-\frac{1}{A_n(s)B_n(s)}\right)
-\frac{s^{-2}\varepsilon_s}{B_n(s)}\\
&-\frac{s^{-2}d_s^2}{A_n(s)B_n(s)}
-\frac{s^{-4}\delta_s}{A_n(s)B_n(s)},
\end{align*}
then
\begin{equation}\label{eq:asymptotic-decomp}
f_n(s)-\frac{s^{-2}}{A_n(s)}
=
\frac{b_4s^{-6}}{A_n(s)^2}
-\frac{b_3^2s^{-8}}{A_n(s)^3}
+\xi_n(s).
\end{equation}

\medskip
\noindent\textbf{Estimate of the remainder.}
From \eqref{eq:eps-bound} and \eqref{eq:eta-bound}, there exists some constant $C_\mu>0$ such that
\[
|d_s|\leq C_\mu s^{-4},
\qquad
|R_s|\leq C_\mu s^{-2}.
\]
Since $A_n(s)\geq s^{-2}$, for $M$ sufficiently large,
\[
|d_s|\leq \frac{A_n(s)}{2},
\qquad
s^{-2}R_s^2\leq \frac{A_n(s)^2}{4}.
\]
Hence
\[
cA_n(s)^2\leq B_n(s)\leq CA_n(s)^2
\]
for some constants $c,C>0$ depending only on $\mu$. Moreover,
\[
B_n(s)-A_n(s)^2
=
2A_n(s)d_s+d_s^2+s^{-2}R_s^2,
\]
so that
\[
|B_n(s)-A_n(s)^2|
\leq C_\mu A_n(s)s^{-4}.
\]
Consequently,
\[
\left|
\frac{1}{B_n(s)}-\frac{1}{A_n(s)^2}
\right|
\leq
C_\mu\frac{s^{-4}}{A_n(s)^3},
\qquad
\left|
\frac{1}{A_n(s)^3}
-\frac{1}{A_n(s)B_n(s)}
\right|
\leq
C_\mu\frac{s^{-4}}{A_n(s)^4}.
\]

For fixed $M$, let
\[
\theta_M:=
\sup_{s\geq M}s^4\varepsilon_s.
\]
By \eqref{eq:eps-bound}, $\theta_M\to0$ as $M\to\infty$. Using the previous estimates and $|d_s|\leq C_\mu s^{-4}$, we obtain
\[
\begin{aligned}
|\xi_n(s)|
\leq C_\mu\Bigg(
\frac{s^{-10}}{A_n(s)^3}
+\frac{s^{-12}}{A_n(s)^4}
+\theta_M\frac{s^{-6}}{A_n(s)^2}+\frac{s^{-10}}{A_n(s)^3}
+\frac{s^{-9}}{A_n(s)^3}
\Bigg).
\end{aligned}
\]
With the change of variable $u=s^{-2}$, we have
\[
\frac{s}{n}\,ds=-\frac{1}{2n}u^{-2}\,du,
\]
and it follows that
\[
\begin{aligned}
\int_M^\infty\frac{s}{n}|\xi_n(s)|\,ds
\leq \frac{C_\mu}{n}\Bigg[
&\int_0^{M^{-2}}\frac{u^3}{(u+n^{-1})^3}\,du+\int_0^{M^{-2}}\frac{u^4}{(u+n^{-1})^4}\,du
+\theta_M\int_0^{M^{-2}}
\frac{u}{(u+n^{-1})^2}\,du\\ &+\int_0^{M^{-2}}\frac{u^3}{(u+n^{-1})^3}\,du+\int_0^{M^{-2}}\frac{u^{5/2}}{(u+n^{-1})^3}\,du
\Bigg].
\end{aligned}
\]
The first, second, fourth and fifth integrals are uniformly bounded in $n$, while
\[
\int_0^{M^{-2}}
\frac{u}{(u+n^{-1})^2}\,du
=
\log n+O_M(1).
\]
Thus
\begin{equation}\label{eq:remainder-o}
\limsup_{n\to\infty}
\frac{n}{\log n}
\int_M^\infty\frac{s}{n}|\xi_n(s)|\,ds
\leq C_\mu\theta_M.
\end{equation}

\medskip
\noindent\textbf{Integration of the leading terms.}
From \eqref{eq:asymptotic-decomp}, the two leading contributions are
\[
\begin{aligned}
T_{4,n}(M)
&:=
\int_M^\infty
\frac{s}{n}\frac{b_4s^{-6}}{A_n(s)^2}\,ds\\
&=
\frac{nb_4}{2}
\int_{M^2}^\infty\frac{dy}{y(n+y)^2},
\end{aligned}
\]
and
\[
\begin{aligned}
T_{3,n}(M)
&:=
-\int_M^\infty
\frac{s}{n}\frac{b_3^2s^{-8}}{A_n(s)^3}\,ds\\
&=
-\frac{n^2b_3^2}{2}
\int_{M^2}^\infty\frac{dy}{y(n+y)^3},
\end{aligned}
\]
where $y=s^2$. Using the elementary decompositions
\[
\frac{1}{y(n+y)^2}
=
\frac{1}{n^2}\frac{1}{y}
-\frac{1}{n^2}\frac{1}{n+y}
-\frac{1}{n}\frac{1}{(n+y)^2},
\]
and
\[
\begin{aligned}
\frac{1}{y(n+y)^3}
=&
\frac{1}{n^3}\frac{1}{y}
-\frac{1}{n^3}\frac{1}{n+y}-\frac{1}{n^2}\frac{1}{(n+y)^2}
-\frac{1}{n}\frac{1}{(n+y)^3},
\end{aligned}
\]
we obtain
\begin{equation}\label{eq:T4}
T_{4,n}(M)
=
\frac{b_4}{2}\frac{\log n}{n}
+O_M(n^{-1}),
\end{equation}
and
\begin{equation}\label{eq:T3}
T_{3,n}(M)
=
-\frac{b_3^2}{2}\frac{\log n}{n}
+O_M(n^{-1}).
\end{equation}

Combining \eqref{eq:compact}, \eqref{eq:remainder-o}, \eqref{eq:T4}, and \eqref{eq:T3}, we deduce that
\[
\limsup_{n\to\infty}\left|\frac{n}{\log n}\Gamma(\mu_n\mid\mathbf b)-\frac{b_4-b_3^2}{2}\right|
\leq C_\mu\theta_M.
\]
Finally, letting $M\to\infty$ and using $\theta_M\to0$, we conclude that
\[
\Gamma(\mu_n\mid\mathbf b) = \frac{b_4-b_3^2}{2}\frac{\log n}{n} +o\left(\frac{\log n}{n}\right),
\]
as claimed.
\end{proof}

\begin{example}[Sharpness of the entropic rate]\label{ex:four-point}
Fix $d\in(0,1)$ and consider
\[
 \nu_d=\frac14\bigl(\delta_{\sqrt{1+d}}+\delta_{-\sqrt{1+d}}+\delta_{\sqrt{1-d}}+\delta_{-\sqrt{1-d}}\bigr).
\]
Write $\nu_{d,n}=\D_{n^{-1/2}}(\nu_d^{\uplus n})$. Then
\[
 \Gamma(\nu_{d,n})=-\frac{d^2}{2}\frac{\log n}{n}+O\!\left(\frac1n\right),
\]
\end{example}

\begin{proof}
A direct calculation gives
\[
 G_{\nu_d}(w)=\frac{w(w^2-1)}{(w^2-1)^2-d^2},
 \qquad
 F_{\nu_d}(w)=w-\frac1w-\frac{d^2}{w(w^2-1)}.
\]
Using the Boolean convolution and dilation identities,
\[
 F_{\nu_{d,n}}(z)
 =z-\frac1z-\frac{d^2}{z(nz^2-1)},
 \]
and therefore
\[
 G_{\nu_{d,n}}(z)
 =\frac{z(nz^2-1)}{(z^2-1)(nz^2-1)-d^2}.
\]
Set
\[
 \Delta_n=\sqrt{(n-1)^2+4nd^2},
 \qquad
 u_{\pm,n}=\frac{n+1\pm\Delta_n}{2n},
\]
and
\[
 \omega_{+,n}=\frac{n-1+\Delta_n}{2\Delta_n},
 \qquad
 \omega_{-,n}=\frac{\Delta_n-(n-1)}{2\Delta_n}.
\]
Partial fractions give
\[
 \nu_{d,n}
 =\frac{\omega_{+,n}}{2}(\delta_{\sqrt{u_{+,n}}}+\delta_{-\sqrt{u_{+,n}}})
 +\frac{\omega_{-,n}}{2}(\delta_{\sqrt{u_{-,n}}}+\delta_{-\sqrt{u_{-,n}}}).
\]
Hence
\[
 \Gamma(\nu_{d,n})
 =\frac{\omega_{+,n}}{2}\log u_{+,n}+\frac{\omega_{-,n}}{2}\log u_{-,n}.
\]
As $n\to\infty$,
\[
 \Delta_n=n-1+2d^2+O(n^{-1}),
\]
\[
 u_{+,n}=1+\frac{d^2}{n}+O(n^{-2}),
 \qquad
 u_{-,n}=\frac{1-d^2}{n}+O(n^{-2}),
\]
and
\[
 \omega_{+,n}=1-\frac{d^2}{n}+O(n^{-2}),
 \qquad
 \omega_{-,n}=\frac{d^2}{n}+O(n^{-2}).
\]
Substitution yields
\[
 \Gamma(\nu_{d,n})=-\frac{d^2}{2}\frac{\log n}{n}+O(n^{-1}).
\]
Since $\nu_{d,n}\in\mathcal{P}_{0,1}(\mathbb{R})$, the identity \eqref{eq:relative entropy} gives
\[
 \Gamma(\nu_{d,n}\mid \mathbf{b})
 =\frac{d^2}{2}\frac{\log n}{n}+O(n^{-1}),
\]
as claimed.
\end{proof}

\subsection{Quantitative Fisher information asymptotics}\label{sec:fisher}

We now prove Theorem~\ref{thm:fisher-asymptotic}.

\begin{proof}
Let $\rho=\mu^{(2)}$. Then $\rho$ is supported on $\mathbb{R}_+$, $\int x\,d\rho(x)=1$, and
\[
 \int(x-1)^2\,d\rho(x)=m_4(\mu)-1.
\]
For any $s>0$, set
\[
 h(s)=\int_0^\infty\frac{d\rho(x)}{s+x},
 \qquad
 r(s)=\frac1{h(s)}-s,
 \qquad
 q(s)=1-r(s).
\]
Jensen's inequality gives
\[
 \frac1{s+1}\le h(s)<\frac1s,
\]
so $0<r(s)\le1$ and $0\le q(s)<1$. By \eqref{eq:square-convolution},
\[
 \mu_n^{(2)}=\D_{1/n}(\rho^{\uplus n}),
 \qquad
 K_{\mu_n^{(2)}}(z)=K_\rho(nz).
\]
Hence, for any $s>0$,
\[
 -G_{\mu_n^{(2)}}(-s)=\frac1{s+r(ns)}.
\]
Applying the change of variables $u = ns$ to the identity $\int_0^\infty (s+1)^{-2} \, ds = 1$, and setting $b_n(u) \coloneqq 1 + u/n$, we have
\begin{equation*}
    1 = \frac{1}{n} \int_0^\infty \frac{du}{b_n(u)^2}.
\end{equation*}
Substituting the identity in the integral representation of $\Psi(\mu_n)$ in Lemma~\ref{lem:stieltjes} yields
\begin{equation}\label{eq:psi-integral}
    \Psi(\mu_n) - 1 = \frac{1}{n} \int_0^\infty H_n(u) \, du,
\end{equation}
where $H_n(u) \coloneqq (b_n(u) - q(u))^{-2} - b_n(u)^{-2}$. A straightforward calculation gives
\begin{align}\label{eq:Hn}
    H_n(u) &= \frac{b_n(u)^2 - (b_n(u) - q(u))^2}{b_n(u)^2 (b_n(u) - q(u))^2} \nonumber \\
    &= q(u) \frac{2b_n(u) - q(u)}{b_n(u)^2 (b_n(u) - q(u))^2},
\end{align}
Recall that $b_4=m_4(\mu)-1.$ We first establish the following asymptotic behavior of $q$:
\begin{equation}\label{eq:q-asymptotic}
    \lim_{u\to \infty}uq(u)=b_4.
\end{equation}
Indeed, a rearrangement gives
\[
 q(u)=\frac{(u+1)h(u)-1}{h(u)},
\]
and by $\int(x-1)d\rho(x)=0$, we have
\[
 uq(u)
 =\frac{\displaystyle\int_0^\infty\frac{u(x-1)^2}{u+x}\,d\rho(x)}{\displaystyle(u+1)\int_0^\infty\frac1{u+x}\,d\rho(x)}.
\]
The dominated convergence theorem yields \eqref{eq:q-asymptotic}. If $b_4=0$, then $\rho=\delta_1$ and symmetry forces $\mu=\mathbf{b}$, so the conclusion is immediate. Assume henceforth $b_4>0$, and it remains to prove
\begin{equation}\label{eq:Hn-asymptotic}
    \frac{1}{n}\int_0^\infty H_n(u)\,du
 =2b_4\frac{\log n}{n}+o\left(\frac{\log n}{n}\right).
\end{equation}

We split the integral into a bounded part, an intermediate bulk, and a tail.

\medskip
\noindent\textbf{Bounded region.}
From \eqref{eq:Hn}, $0\le H_n(u)\le r(u)^{-2}.$ Moreover, $r(s)/s=(sh(s))^{-1}-1$ is decreasing in $s$, since $sh(s)$ is increasing in $s$. Thus for any $0<u\le M$,
\[
 r(u)\ge\frac{r(M)}M u,
 \qquad
 u+r(u)\le\left(1+\frac M{r(M)}\right)r(u).
\]
Since
\[
 \Psi(\mu)=\int_0^\infty\frac{du}{(u+r(u))^2}<\infty,
\]
we obtain
\begin{equation}\label{eq:bounded-Hn}
    \sup_{n\ge1}\int_0^M H_n(u)\,du<\infty.
\end{equation}

\medskip
\noindent\textbf{Bulk region.}
Fix $\eta\in(0,b_4)$. By \eqref{eq:q-asymptotic}, there exists some sufficiently large $M$ such that
\[
 \frac{b_4-\eta}{u}\le q(u)\le\frac{b_4+\eta}{u}
\qquad (u\ge M).
\]
For $M\le u\le\delta n$, write
\[
 H_n(u)=q(u)\Phi(b_n(u),q(u)),
 \qquad
 \Phi(b,q)=\frac{2b-q}{b^2(b-q)^2}.
\]
Because $\Phi(1,0)=2$, one may first choose $M$ large and then $\delta>0$ small enough that
\[
 2-\eta\le \Phi(b_n(u),q(u))\le2+\eta
\]
uniformly on $[M,\delta n]$. Hence
\begin{equation}\label{eq:bulk-Hn}
    (2-\eta)(b_4-\eta)
 \log\frac{\delta n}{M}
 \le\int_M^{\delta n}H_n(u)\,du
 \le(2+\eta)(b_4+\eta)
 \log\frac{\delta n}{M}.
\end{equation}

\medskip
\noindent\textbf{Tail.}
Again by \eqref{eq:q-asymptotic}, there exists $C>0$ such that $q(u)\le C/u$ for $u\ge M$. For fixed $\delta>0$ and sufficiently large $n$, $u\ge\delta n$ implies $q(u)\le1/2$, and hence
\[
 H_n(u)\le\frac{8q(u)}{(1+u/n)^3}
 \le\frac{8C}{u(1+u/n)^3}.
\]
Therefore
\begin{equation}\label{eq:tail-Hn}
    \int_{\delta n}^\infty H_n(u)\,du
 \le 8C\int_\delta^\infty\frac{dv}{v(1+v)^3}
 =O_\delta(1).
\end{equation}
Dividing \eqref{eq:bounded-Hn}, \eqref{eq:bulk-Hn}, and \eqref{eq:tail-Hn} by $ n$, then letting $n\to\infty$ and $\eta\downarrow0$, yields \eqref{eq:Hn-asymptotic}. Substitution into \eqref{eq:psi-integral} gives
\[
 \Psi(\mu_n)-1
 =2b_4\frac{\log n}{n}+o\!\left(\frac{\log n}{n}\right),
\]
as claimed.
\end{proof}

\begin{example}[Sharpness of the Boolean Fisher information rate]\label{ex:fisher-four-point}
Under the setup and notation of Example~\ref{ex:four-point}, the Boolean Fisher information of $\nu_{d,n}$ satisfies
\[
\Psi(\nu_{d,n}) - 1 = 2d^2 \frac{\log n}{n} + O\!\left(\frac1n\right).
\]
Using the explicit discrete representation of $\nu_{d,n}$ from Example~\ref{ex:four-point}, the Boolean Fisher information evaluates directly to
\[
\Psi(\nu_{d,n}) = \omega_{+,n}^2 \frac{1}{u_{+,n}} + 2 \omega_{+,n} \omega_{-,n} \frac{\log u_{+,n} - \log u_{-,n}}{u_{+,n} - u_{-,n}} + \omega_{-,n}^2 \frac{1}{u_{-,n}}.
\]
Substituting the asymptotic expansions of $u_{\pm,n}$, $\omega_{+,n}$, and $\omega_{-,n}$ as $n \to \infty$, we expand each term separately:
\begin{align*}
   \omega_{+,n}^2 \frac{1}{u_{+,n}} = \left(1 - \frac{2d^2}{n} + O(n^{-2})\right) \left(1 - \frac{d^2}{n} + O(n^{-2})\right) = 1 - \frac{3d^2}{n} + O(n^{-2}).
\end{align*}
\begin{align*}
   2 \omega_{+,n} \omega_{-,n} \frac{\log u_{+,n} - \log u_{-,n}}{u_{+,n} - u_{-,n}} 
   &= \left(\frac{2d^2}{n} + O(n^{-2})\right) \frac{\log n - \log(1-d^2) + O(n^{-1})}{1 + O(n^{-1})} \\
   &= 2d^2 \frac{\log n}{n} - \frac{2d^2 \log(1-d^2)}{n} + O\!\left(\frac{\log n}{n^2}\right).
\end{align*}
\begin{align*}
    \omega_{-,n}^2 \frac{1}{u_{-,n}} 
   &= \left(\frac{d^4}{n^2} + O(n^{-3})\right) \left(\frac{n}{1-d^2} + O(1)\right) = \frac{d^4}{(1-d^2)n} + O(n^{-2}).
\end{align*}
Summing these three contributions gives
\[
\Psi(\nu_{d,n}) = 1 + 2d^2 \frac{\log n}{n} + O\!\left(\frac1n\right).
\]
Using the identity \eqref{eq:relative fisher} yields the desired result.
\end{example}

\bibliographystyle{plain}
\bibliography{references}
\end{document}